\documentclass[a4paper,12pt,oneside]{amsart}

\usepackage[all]{xy}
\usepackage{amsfonts}
\usepackage{amssymb}
\usepackage{amsmath}
\usepackage{latexsym}
\input cyracc.def
\font\tencyr=wncyr10
\def\cyr{\tencyr\cyracc}
\newcommand{\ts}{\mbox{\cyr Sh}}
\newcommand{\ov}{\overline}
\newcommand{\Q}{\mathbb{Q}}
\newcommand{\Z}{\mathbb{Z}}
\newcommand{\F}{\mathcal{F}}
\newcommand{\X}{\mathcal{X}}
\renewcommand{\S}{\mathcal{S}}
\newcommand{\Ql}{{\mathbb{Q}_l}}
\newcommand{\Zl}{{\mathbb{Z}_l}}
\newcommand{\Zp}{{\mathbb{Z}_p}}
\newcommand{\Qp}{{\mathbb{Q}_p}}
\newcommand{\N}{\mathbb{N}}
\newcommand{\Zpn}{{\mathbb{Z}_p^{\mathbb{N}}}}
\newcommand{\ol}{\mathcal O}
\newcommand{\m}{\mathfrak m}

\newcommand{\Fv}{{\mathbb{F}_v}}
\newcommand{\il}[1]{\lim_{\buildrel \longleftarrow\over{#1}}}
\newcommand{\dl}[1]{\lim_{\buildrel \longrightarrow\over{#1}}}
\newcommand{\sri}{\twoheadrightarrow}
\newcommand{\iri}{\hookrightarrow}
\newcommand{\ri}{\rightarrow}

\newcommand{\p}{\mathfrak{p}}
\renewcommand{\L}{\Lambda}
\renewcommand{\t}{\tilde}
\newcommand{\wt}{\widetilde}

\newcommand{\g}{\gamma}
\newcommand{\G}{\Gamma}

\renewcommand{\k}{\kappa}
\newtheorem{defin}{Definition}[section]
\newtheorem{prop}[defin]{Proposition}
\newtheorem{conj}[defin]{Conjecture}
\newtheorem{lemma}[defin]{Lemma}
\newtheorem{cor}[defin]{Corollary}
\newtheorem{rem}[defin]{Remark}
\newtheorem{thm}[defin]{Theorem}
\newcommand{\liminv}{\displaystyle \lim_{\leftarrow}}
\newcommand{\limdir}{\displaystyle \lim_{\rightarrow}}

\begin{document}

\title[Control theorems ...]{Control theorems for elliptic curves over function fields}

\author{A. Bandini, I. Longhi}

\begin{abstract} 
Let $F$ be a global field of characteristic $p>0$,
$\F/F$ a Galois extension with $Gal(\F/F)\simeq \Z_p^\N$ and $E/F$ a 
non-isotrivial elliptic curve. We study the behaviour of Selmer groups $Sel_E(L)_l$
($l$ any prime) as $L$ varies through the subextensions of $\F$ via appropriate
versions of Mazur's Control Theorem. In the case $l=p$ we let $\F=\bigcup \F_d$ where
$\F_d/F$ is a $\Z_p^d$-extension. With a mild hypothesis on $Sel_E(F)_p$ (essentially 
a consequence of the Birch and Swinnerton-Dyer conjecture) we prove that $Sel_E(\F_d)_p$ is a 
cofinitely generated $\Zp[[Gal(\F_d/F)]]$-module 
and we associate to its Pontrjagin dual a Fitting ideal. This allows to define an
algebraic $L$-function associated to $E$ in $\Zp[[Gal(\F/F)]]$, providing an ingredient for
a function field analogue of Iwasawa's Main Conjecture for elliptic curves.
\end{abstract}

\maketitle

\section{Introduction}

\subsection{Some motivation} Iwasawa theory for elliptic curves over number fields is by now an established subject, 
well developed both in its analytic and algebraic sides. By contrast, its analogue over global function fields is still 
at its beginnings: as far as the authors know, up to now only analytic aspects have been investigated. Our goal in this 
paper is to provide some first steps into understanding the algebraic side as well. To relate the two points of view
we end by proposing a first version of Iwasawa's Main Conjecture, which we hope to investigate
in some future work. We are conscious that our Conjecture \ref{MC} is too coarsely formulated to be completely 
satisfactory; however, our point is mainly to show that it is possible to ask questions of this kind also in 
the characteristic $p$ setting.

We fix $F$ a function field of transcendence degree 1 over its constant field $\mathbb{F}_q\,$, $q$ a power of a
prime $p$, and an elliptic curve $E/F$; we assume that $E$ is non-isotrivial (i.e., $j(E)\notin\mathbb F_q$). In
particular, $E$ has bad reduction at some place of $F$; replacing, if needed, $F$ by a finite extension, we can
(and will) assume that $E$ has good or split multiplicative reduction at any place $v$ of $F$.

\subsection{Analytic theory} We briefly review the state of the art.

\subsubsection{The extensions}\label{cycext}
Let $\wt\F/F$  be an infinite Galois 
extension such that $Gal(\wt\F/F)$ contains a finite index 
subgroup $\G$ isomorphic to $\Zpn$ (an infinite product of $\Zp$'s). 
The reader is reminded that class field theory provides lots of such 
extensions, thanks to the fact that if $L$ is a local field in 
characteristic $p$ and $U_1(L)$ denotes its group of $1$-units then 
$U_1(L)\simeq\Zpn$ (see e.g. \cite[II.5.7]{Ne}). Observe that, exactly 
for this reason, in the function field setting it becomes quite 
natural to concentrate on such a $\G$ rather than on a finite dimensional 
$p$-adic Lie group.

A good example, which closely parallels the classical cyclotomic $\Zp$-extension of a number field, is the 
``cyclotomic extension at $\p$''. In the simplest formulation, we take $F:=\mathbb F_q(T)$, $A:=\mathbb F_q[T]$
and let $\Phi$ be the Carlitz module (see e.g. \cite[Chapter 12]{Ro}). Choose $\p$ a prime of $A$ and, 
for any positive integer $n$, let $\Phi[\p^n]$ denote the $\p^n$-torsion of the Carlitz module. Let $F(\Phi[\p^n])$ be 
the extension of $F$ obtained via the $\p^n$-torsion and let
$\wt\F:=F(\Phi[\p^\infty])=\bigcup F(\Phi[\p^n])\,$.
It is well known that $F(\Phi[\p^n])/F$ is a Galois extension and that
\[ Gal( F(\Phi[\p^n])/F)\simeq (A/\p^n)^*\ ,\]
\[ Gal(\wt\F/F)\simeq \il{n} (A/\p^n)^*
\simeq \Zpn \times (A/\p)^*\ .\]

\subsubsection{The ``$p$-adic $L$-function''}\label{pLfunc} Here by $p$-adic $L$-function 
we mean an element in the Iwasawa algebra $\Z_p[[ Gal(\wt\F/F)]]$ 
(identified with the algebra of $\Zp$-valued measures on $Gal(\wt\F/F)$): 
to our knowledge, up to now no satisfactory closer analogue of the usual 
$p$-adic $L$-function arising in characteristic $0$ was found. 
(A key problem seems to be the lack of an adequate theory of Mellin 
transform, in spite of some attempts by Goss - see \cite{Go2}.)

The first instance of construction of a measure related to $E$ is due to Teitelbaum (\cite[pag. 290-292]{Te}):
the $\wt\F$ he implicitly considers is exactly the Carlitz cyclotomic extension at $\p$ described 
above (where $\p$ is a prime of split multiplicative reduction for $E$). 
Other examples (which can be loosely described as cyclotomic and anticyclotomic at $\infty$) 
were given in \cite{Lo}. For a more detailed discussion see section \ref{pLf}.

\subsection{The present work} Since in this paper we are not going to work out a 
comparison with the analytic theory, our attention 
will be focussed only on a $\Zpn$-extension  $\F/F$, i.e.  
a Galois extension such that $\G:=Gal(\F/F)\simeq\Zpn\,$. 
For example, in the situation described above, one can take $\F$ to be the subfield of  
$F(\Phi[\p^\infty])$ fixed by $(A/\p)^*\,$: then $\F/F$ is a $\Zpn$-extension. We shall consider  
$\Z_p^d$-extensions $\F_d/F$ such that $\F=\bigcup \F_d\,$.

Denote by $\Lambda:=\Zp[[\G]]$ and by $\L_d:=\Zp[[Gal(\F_d/F)]]$ the associated Iwasawa algebras.

In section \ref{SelGr} we will define the ($p$-part of the) Selmer group 
$Sel_E(L)_p$, $L$ any algebraic extension of $F$. For any $d$, let $\S_d$ (resp. $\S$) be the 
Pontrjagin dual of $Sel_E(\F_d)_p$ (resp. of $Sel_E(\F)_p$): it is a $\L_d$-module (resp. a $\L$-module). 
The main result of this paper is the following.

\begin{thm}\label{ThmIntr} 
Assume that all ramified primes in $\F_d/F$ are of split multiplicative reduction for $E$ and
that $Sel_E(F)_p$ is cofinitely generated as a $\Zp$-module. Then $\S_d$ is a finitely generated $\L_d$-module. 
\end{thm}

\begin{rem} \emph{The assumption on $Sel_E(F)_p$ is a consequence of the Birch and Swinnerton-Dyer conjecture. 
Some evidence for the latter can be found in \cite{Brw} and \cite{Ul2}.}
\end{rem}

We get Theorem \ref{ThmIntr} as Corollary \ref{ControlpCor} of our Theorem \ref{Controlp} 
which is an analogue of Mazur's classical Control Theorem (for which the reader is
referred to \cite{Ma} or \cite{Gr2}). \\ 
The most interesting consequence is that it is possible to define a kind of 
Fitting ideal $\wt{Fitt}_\L(\S)\subset\L$ of $\S$ as the intersection of counterimages 
of the $Fitt_{\L_d}(\S_d)$'s in $\L$ (see section \ref{LimFitt}, Definition \ref{DefFittId}).
The natural next step is the formulation of a Main Conjecture relating
this $\wt{Fitt}_\L(\S)$ and a ``$p$-adic $L$-function'' as described above (Conjecture \ref{MC}).

By the same techniques, and with less effort, we also investigate the variation of the $l$-part of the Selmer group
($l$ a prime different from $p$) in subextensions of $\F/F$. Theorem \ref{Controll} shows that also in this case
we can control the Selmer groups; however we lack a good theory of modules over the ring
$\Z_l[[\G]]$ and thus are unable to say much more.

\begin{rem}\label{supsing}
{\em It might be worthwhile to remark that the prime number $p$ is not a 
place of our field, contrary to the classical situation. In particular we don't have to 
ask anything about the reduction of $E$ being not supersingular at some place.}
\end{rem}

\subsection{Structure of this paper} In section \ref{SelGr} we establish notations and define 
our Selmer groups by flat cohomology, which, for $l\neq p$, reduces to the usual Galois cohomology.
Section \ref{lcontrol} is dedicated to the easier case $l\neq p$: here we can establish the 
control theorem without any assumption on the $\Zpn$-extension $\F/F$.
On the contrary the control theorems for $l=p$ are proven (in section \ref{pcontrol}) only for a 
$\Z_p^d$-extension $\F_d/F$: the
technical reasons for this limitation are explained in Remark \ref{kerbnrem}. We begin with Theorem
\ref{Controlp} which regards classical Selmer groups and leads to finitely generated $\L_d$-modules
in Corollary \ref{ControlpCor}. To provide some examples of torsion modules we also
present another control theorem for some slightly modified Selmer groups (Theorem \ref{ControlModSel})
where we change the local conditions at the ramified primes.
In section \ref{algL} we go back to $\F/F$. We define 
the $\L$-ideal $\wt{Fitt}_\L(\S)$, which provides our candidate for an algebraic $L$-function and by 
which we formulate our version of Iwasawa's Main Conjecture.

\noindent {\bf Acknowledgements.} We thank Massimo Bertolini, Henri Darmon, Ralph Greenberg, Marc-Hubert Nicole,
Tadashi Ochiai, David Vauclair and Stefano Vigni for useful suggestions and comments on early drafts of this paper.\\
The first author has been supported by a scientific grant from the 
University of Pisa. The second author has been supported by a scientific grant from the University of Milano.

\section{Selmer groups} \label{SelGr}

\subsection{Notations} For the convenience of the reader we list some
notations that will be used in this paper and describe the basic setting
in which the theory will be developed.

\subsubsection{Fields} Let $L$ be a field: then $\ov{L}$ will denote an algebraic closure and $L^{sep}\subset
\ov{L}$ a separable closure; besides, we put $G_L:=Gal(L^{sep}/L)$.\\ 
If $L$ is a global field (or an algebraic extension of such), $\mathcal M_L$ will be 
its set of places. For any place $v\in\mathcal M_L$ we let
$L_v$ be the completion of $L$ at $v$, $\ol_v$ the ring of integers of 
$L_v\,$, $\m_v$ the maximal ideal of $\ol_v$
and $\mathbb L_v:=\ol_v/\m_v$ the residue field. Let $ord_v$ be the valuation associated to $v$.\\ If $L$ is a
local field, the notations above will often be changed to $\ol_L$, $\m_L$, $\mathbb L$; besides
$U_1(L)\subset\ol_L^{\ast}$ will be the group of $1$-units.

\noindent As stated in the introduction, we fix a global field $F$ of characteristic $p>0$ and
an algebraic closure $\ov{F}$. For any place
$v\in\mathcal M_F$ we choose $\ov{F_v}$ and an embedding $\ov{F}\hookrightarrow\ov{F_v}$, so to get a restriction map
$G_{F_v}\hookrightarrow G_F$. All algebraic extensions of $F$ (resp. $F_v$) will be assumed to
be contained in $\ov{F}$ (resp. $\ov{F_v}$).

Script letters will denote infinite extensions of $F$: more precisely we fix a $\Zpn$-extension $\F/F$ 
with Galois group $\G$ and for any $d\ge 1$ we will consider $\F_d \subset \F$, a $\mathbb{Z}_p^d$-extension of $F$.

\subsubsection{Elliptic curves} We fix an elliptic curve $E/F$, non-isotrivial and having split multiplicative
reduction at all places supporting its conductor. The reader is reminded that then at such 
places $E$ is isomorphic to a Tate curve, i.e. $E(F_v)\simeq F_v^{\ast}/q_{E,v}^{\Z}$ for 
some $q_{E,v}$.\\
For any positive integer $n$ let
$E[n]$ be the scheme of $n$-torsion points. Moreover, for any prime $l$, let
$E[l^\infty]:=\limdir E[l^n]\,$.\\
For any $v\in\mathcal{M}_F$ we choose a minimal Weierstrass equation for $E$.
Let $E_v$ be the reduction of $E$ modulo $v$ and for any point $P\in E$ let
$P_v$ be its image in $E_v\,$.\\
Besides $E_{v,ns}(\Fv)$ is the set of nonsingular points of $E_v(\Fv)$, and
\[ E_0(F_v) := \{\,P\in E(F_v)\,\mid\ P_v\in E_{v,ns}(\Fv)\,\}\ . \]
By the theory of Tate curves we know that, in case of bad reduction, $E_0(F_v)\simeq \ol_v^*$ (see e.g. \cite[V section 4]{Si2}).\\
Finally, let $T_v:=E(F_v)/E_0(F_v)\,$. In case of bad reduction $T_v$ is the group of components 
of the special fibre and our hypothesis implies that its order is $-ord_v(j(E))$ (see e.g. \cite[IV.9.2]{Si2}). \\
For all basic facts about elliptic curves, the reader is referred to Silverman's books \cite{Si1} and \cite{Si2}.

\subsubsection{Duals} For $X$ a topological abelian group, we denote its Pontrjagin dual by $X^{\vee}:=Hom_{cont}(X,\mathbb{C}^{\ast})$. 
In the cases considered in this paper, $X$ will be a (mostly discrete) topological $\Zl$-module, so that $X^{\vee}$ 
can be identified with $Hom_{cont}(X,\Ql/\Zl)$ and it has a natural structure of $\Zl$-module.\\
The reader is reminded that to say that an $R$-module $X$ ($R$ any ring) is cofinitely generated means that $X^{\vee}$ is a 
finitely generated $R$-module.

\subsection{The Selmer groups}\label{Selmer}
We are interested in torsion subschemes of the elliptic curve $E$.
Since $char\,F=p$, in order to deal with the $p$-torsion and to define Selmer groups with the usual cohomological
techniques, we need to consider flat cohomology of group schemes.\\
For the basic theory of sites and cohomology on a site see \cite[Chapters II, III]{Mi1}. Briefly, for any scheme
$X$ we let $X_{fl}$ be the subcategory of ${\bf Sch}/X$ (schemes over $X$) whose structure morphisms are locally of
finite type. Moreover $X_{fl}$ is endowed with the flat topology, i.e., if we let $Y\ri X$ be an element of
$X_{fl}\,$, then a covering of $Y$ is a family $\{g_i:U_i\ri Y\,\}$ such that $Y=\bigcup g_i(U_i)$ and each
$g_i$ is a flat morphism locally of finite type.\\ 
We only consider flat cohomology so when we write a
scheme $X$ we always mean $X_{fl}\,$. 

\begin{defin}\label{defflco}
Let $\mathcal{P}$ be a sheaf on $X$ and consider the global section functor
sending $\mathcal{P}$ to $\mathcal{P}(X)$.
The $i$-th flat cohomology group of $X$ with values in $\mathcal{P}$, denoted by $H_{fl}^i(X,\mathcal{P})\,$, is the
value at $\mathcal{P}$ of the $i$-th right derived functor of the global section functor.
\end{defin}

Let $L$ be an algebraic extension of $F$ and $X_L:=Spec\,L$. For any positive integer $m$ the
group schemes $E[m]$ and $E$ define sheaves on $X_L$ (see \cite[II.1.7]{Mi1}): for example 
$E[m](X_L):=E[m](L)$. Consider the exact sequence
\[ E[m]\iri E {\buildrel m\over\longrightarrow\!\!\!\!\!\rightarrow} E \]
and take flat cohomology with respect to $X_L$ to get
\[ E(L)/mE(L)\iri H_{fl}^1(X_L,E[m])\ri H_{fl}^1(X_L,E)\ . \]
In particular let $m$ run through the powers $l^n$ of a prime $l$. Taking direct limits
one gets an injective map (a ``Kummer homomorphism'')
\[ \k : E(L)\otimes \Q_l/\Z_l \iri \dl{n} H_{fl}^1(X_L,E[l^n])=: H_{fl}^1(X_L,E[l^\infty])\ . \]
Exactly as above one can build local Kummer maps for any place 
$v\in\mathcal{M}_L$
\[ \k_v : E(L_v)\otimes \Q_l/\Z_l \iri H_{fl}^1(X_{L_v},E[l^\infty]) \]
where $X_{L_v}:=Spec\,L_v\,$.

\begin{defin}\label{Sel}
The \emph{$l$-part of the Selmer group} of $E$ over $L$, denoted by
$Sel_E(L)_l\,$, is defined to be
\[ Sel_E(L)_l:=Ker\{ H_{fl}^1(X_L,E[l^\infty])\ri\prod_{v\in\mathcal{M}_L}
H_{fl}^1(X_{L_v},E[l^\infty])/Im\,\k_v \,\} \]
where the map is the product of the natural restrictions between cohomology groups.\\ 
\end{defin}   

The reader is reminded that the \emph{Tate-Shafarevich group} 
$\ts(E/L)$ fits into the exact sequence
\[ E(L)\otimes \Q_l/\Z_l \iri Sel_E(L)_l \sri \ts(E/L)[l^\infty]\ .\]
According to the function field version of the Birch and 
Swinnerton-Dyer conjecture, $\ts(E/L)$ is finite for any finite extension
$L$ of $F$. Applying to this last 
sequence the exact functor $Hom(\cdot,\Q_l/\Z_l)$, it follows that
\[ rank_{\Z_l}(Sel_E(L)_l^{\vee})=rank_{\Z}(E(L)) \]
(recall that the cohomology groups $H_{fl}^i\,$, hence the Selmer groups, are endowed with the discrete topology).  

Letting $L$ vary through subextensions of $\F/F$, the groups $Sel_E(L)_l$ admit natural actions by $\mathbb{Z}_l$ (because of
$E[l^\infty]\,$) and by $\G$. Hence they are modules over the Iwasawa algebra $\mathbb{Z}_l[[\G]]$.

\subsubsection{The case $l=p$}\label{Selp}
We shall consider $\Z_p^d$-extensions $\F_d/F$ such that $\F=\bigcup \F_d\,$.
One of the main reasons for doing this is the fact that $Sel_E(\F_d)_p$ is a $\L_d:=\Zp[[Gal(\F_d/F)]]$-module
and, thanks to the (noncanonical) isomorphism $\L_d\simeq \Zp[[T_1,\dots ,T_d]]$, we have a fairly precise
description of such modules (see, for example, \cite[VII.4]{Bo}). The first step towards the
proof of Theorem \ref{ThmIntr} will be the study of the restriction maps
\[ Sel_E(L)_p\longrightarrow Sel_E(\F_d)_p^{Gal(\F_d/L)} \]
as $L$ varies in $\F_d/F$.

\subsubsection{The case $l\neq p$}\label{Sell}
Here we deal with the full $\Zpn$-extension and again study the maps
\[ Sel_E(L)_l\longrightarrow Sel_E(\F)_l^{Gal(\F/L)} \ .\]\\
Note that to define $Sel_E(L)_l$ we could have used the sequence
\[ E[l^n](\ov{F})\iri E(F^{sep}) {\buildrel l^n\over\longrightarrow\!\!\!\!\!\rightarrow} E(F^{sep}) \]
and classical Galois (=\'etale) cohomology since, in this case,
\[ H_{fl}^1(X_L,E[l^n])\simeq H_{et}^1(X_L,E[l^n])\simeq H^1(G_L,E[l^n](\ov{F})) \]
(see \cite[III.3.9]{Mi1}). In order to lighten notations, each time we work with
$l\neq p$ (i.e., in the next chapter) we shall use the classical notation $H^i(L,\cdot)$ instead of
$H^i(G_L,\cdot)\simeq H^i_{fl}(X_L,\cdot)$ and write $E[n]$ for $E[n](\ov{F})$, putting 
$E[l^{\infty}]:=\bigcup E[l^n]$.\\
In this case the Kummer map
\[ \k : E(L)\otimes \Ql/\Zl \iri H^1(G_L,E[l^\infty]) \]
has an explicit description as follows. Let $\alpha\in E(L)\otimes \Ql/\Zl$ be represented as $\alpha= P\otimes
\frac{a}{l^k}\,$ ($a\in\Z$) and let $Q\in E(L^{sep})$ be such that $aP=l^kQ$. Then $\k(\alpha)=\varphi_\alpha\,$,
where $\varphi_\alpha(\sigma):=\sigma(Q)-Q$ for any $\sigma\in G_L\,$.\\

\section{Control theorem for $l\neq p$} \label{lcontrol}

\subsection{The image of $\k_v$} We start by giving a more precise 
description of $Im\,\k_v$ (following the path
traced by Greenberg in \cite{Gr1} and \cite{Gr2}).

\begin{lemma}\label{funcfieldLutz}
If $E$ has split multiplicative reduction at $v$ then $E(F_v)$ contains a finite index subgroup
isomorphic to $U_1(F_v)$.
\end{lemma}

\noindent\emph{Proof.} The hypothesis implies that $E$ is isomorphic to a Tate curve: in 
particular, $E(F_v)\simeq F_v^{\ast}/q_{E,v}^{\Z}$ for some $q_{E,v}$ (the Tate period). 
Now it suffices to observe that $U_1(F_v)$ embeds into $F_v^{\ast}/q_{E,v}^{\Z}$ with finite 
index $ord_v(q_{E,v})\cdot |\mathbb{F}_v^{\ast}|$.$\qquad\square$

\begin{rem} \emph{The previous lemma can be seen as a (partial) function field analog of
Lutz's Theorem (\cite[VII.6.3]{Si1}). In the characteristic 0 case for a finite extension $K/\Qp$ and an
elliptic curve $E$ defined over $\Qp$ one finds that $E(K)$ contains a subgroup isomorphic to the ring of integers
of $K$ (taken as an additive group), i.e. a subgroup isomorphic to $\mathbb{Z}_p^{[K:\,\Qp]}\,$.\\ 
In the function field (characteristic $p$) case we lack the logarithmic function and we have to deal with the multiplicative
group of the ring $\ol_v\,$. For any complete local field like $F_v$ one has $U_1(F_v)\simeq \Zpn\,$, hence, for
any finite extension $K_w/F_v\,$, one finds both in $E(K_w)$ and in $E(F_v)$ a subgroup of finite index isomorphic to
$\Zpn\,$ making it hard to emphasize any kind of relation with the degree $[K_w:F_v]$.}
\end{rem}

\begin{prop}\label{Imkvl}
Let $l\neq p$: then for all places $v$ of $F$ $$E(F_v)\otimes \Ql/\Zl =0$$ (i.e. $Im\,\k_v =0$).
\end{prop}

\noindent\emph{Proof.} We consider two cases, according to the behaviour 
of $E$ at $v$.

{\bf Case 1:} $E$ has good reduction at $v$.\\ 
In this case $E[l^\infty]\simeq E_v[l^\infty]$ and by the
N\'eron-Ogg-Shafarevich criterion this is an isomorphism of $G_{F_v}$-modules. 
Hence one has a sequence of maps
\[ H^1(F_v,E[l^\infty]){\buildrel \sim\over\longrightarrow} 
H^1(F_v,E_v[l^\infty])\iri H^1(F_v,E_v(\ov{\mathbb{F}_v})) \] 
where the last map on the right is induced by the natural inclusion $E_v[l^\infty]\iri E_v(\ov{\mathbb{F}_v})$. 
Such map is injective because $E_v(\ov{\mathbb{F}_v})$ is a torsion abelian group and $E_v[l^\infty]$ is its $l$-primary part. 
Composing  $\k_v$ with this sequence, one gets an injection 
\[ E(F_v)\otimes \Ql/\Zl \iri H^1(F_v,E_v(\ov{\mathbb{F}_v}))\ .\] 
By definition $\varphi\in Im\,\k_v$ implies $\varphi(\sigma)=\sigma (Q)-Q$ for some 
$Q\in E(F_v^{sep})$, hence the image of $\varphi$ in $H^1(F_v,E_v(\ov{\mathbb{F}_v}))$ is $\varphi_v$ 
with $\varphi_v(\sigma)= \sigma (Q_v)-Q_v$ ($Q_v\in E_v(\ov{\mathbb{F}_v})$ is the reduction of $Q$ mod $v$). 
Thus $\varphi_v$ is a 1-coboundary; since all
the maps involved are injective, one gets $\varphi =0$ and finally $E(F_v)\otimes \Ql/\Zl =0$.

{\bf Case 2:} $E$ has bad reduction at $v$, i.e. split multiplicative (according to our hypothesis).\\
By Lemma \ref{funcfieldLutz} one has
\[ 0\ri U_1(F_v) \ri E(F_v) \ri A \ri 0  \]
where $A$ is a finite group. Since $U_1(F_v)\simeq \Zpn$ one gets
\[ 0=U_1(F_v)\otimes \Ql/\Zl \ri E(F_v)\otimes \Ql/\Zl \ri
A\otimes \Ql/\Zl=0\ .\qquad\square \]

\subsection{The theorem} 
Recall our $\Zpn$-extension $\F/F$: we choose any sequence of
finite extensions of $F$ such that
\[ F=F_0\subset F_1\subset\cdots\subset F_n\subset\cdots\subset
\bigcup F_n=\F\ .\]
In this setting we let $\G_n:=Gal(\F/F_n)$ and $\G=Gal(\F/F)$ so that
\[ \Zpn\simeq \G \simeq \il{n} \G/\G_n = \il{n} Gal(F_n/F) \ .\]

\noindent For example if $\F$ is the cyclotomic extension at $\p$ (see section \ref{cycext}),
one can take $F_n$ to be the subfield of $F(\Phi[\p^n])$ fixed by $(A/\p)^*\,$.

We will denote by $v_n$ (resp. $w$) primes of $F_n$ (resp. $\F$) and, to shorten notations, we shall denote by
$F_{v_n}$ (resp. $\F_{w}$) the completion of $F_n$ at $v_n$ (resp. of $\F$ at $w$). Finally, for any 
algebraic extension $L/F$, we put
$$\mathcal{G}(L):=Im\big\{H^1(L,E[l^\infty])\rightarrow\prod_{v\in\mathcal{M}_L}\,H^1(L_v,E[l^\infty])/Im\,\k_v\big\}\,
.$$

For any $n$, the restriction $Sel_E(F_n)_l\ri Sel_E(\F)_l^{\G_n}$
fits into the following diagram (with exact rows)
\[ \xymatrix{
Sel_E(F_n)_l \ar[d]^{a_n} \ar@{^{(}->}[r] &
H^1(F_n,E[l^\infty]) \ar[d]^{b_n} \ar@{->>}[r] &
\mathcal{G}(F_n)\ar[d]^{c_n}\\
Sel_E(\F)_l^{\G_n} \ar@{^{(}->}[r] &
H^1(\F,E[l^\infty])^{\G_n} \ar[r] &
\mathcal{G}(\F) } \]

\begin{thm}\label{Controll}
For any $n$, $Sel_E(F_n)_l\simeq Sel_E(\F)_l^{\G_n}\,$.
\end{thm}

\noindent\emph{Proof.} The snake lemma applied to the diagram above shows that it is enough to prove that
$Ker\,b_n=Coker\,b_n=Ker\,c_n=0$ (i.e., to prove that $a_n$ is an isomorphism we shall prove that $c_n$ is
injective and that $b_n$ is an isomorphism as well).\\ \emph{The map $b_n\,$.} By the Hochschild-Serre spectral
sequence (see \cite[2.1.5]{NSW}) one has
\[ Ker\,b_n \simeq H^1(\G_n,E[l^\infty]^{G_{\F}}) \ ,\]
and
\[ Coker\,b_n \subseteq H^2(\G_n,E[l^\infty]^{G_{\F}}) \ .\]
The group $E[l^\infty]^{G_{\F}}=E[l^\infty](\F)$ is $l$-primary. Moreover
\[ \G_n=Gal(\F/F_n)\simeq \il{m} Gal(F_m/F_n) \]
is a pro-$p$-group for all $n\ge 0$ and
\[ H^i(\G_n,E[l^\infty](\F))\simeq
\dl{m} H^i(Gal(F_m/F_n),E[l^\infty](F_m))\ .\]
For any $m\ge n\ge 0$, $H^i(Gal(F_m/F_n),E[l^\infty](F_m))$ ($i\ge 1$)
is killed by $|E[l^\infty](F_m)|$, which is a (finite) power of $l$,
and by $[F_m:F_n]$ (by the $cor\circ res$ map, see e.g. \cite[1.6.1]{NSW}) 
which is a power of $p\neq l$. Hence, for $i\ge 1$, 
$H^i(Gal(F_m/F_n),E[l^\infty](F_m))=0$ and eventually
$H^i(\G_n,E[l^\infty](\F))=0$ as well.\\
\emph{The map $c_n\,$.} By Proposition \ref{Imkvl} $Ker\,c_n$ is contained  
in the kernel of the natural map
\[ d_n\,:\, \prod_{v_n\in\mathcal{M}_{F_n}} H^1(F_{v_n},E[l^\infty])\ri
\prod_{w\in\mathcal{M}_\F} H^1(\F_w,E[l^\infty])\ .\]
Considering every component, we find maps
\[ d_w\,:\, H^1(F_{v_n},E[l^\infty])\ri H^1(\F_w,E[l^\infty]) \] 
and $Ker\,c_n\subseteq\prod_{w\in\mathcal{M}_\F} Ker\,d_w\,$.\\ 
As we have seen for $Ker\,b_n\,$, one has from the inflation restriction sequence
\[ Ker\,d_w\simeq H^1(Gal(\F_w/F_{v_n}),E[l^\infty]^{G_{\F_w}}) \] 
i.e., $Ker\,d_w= 0$ because $Gal(\F_w/F_{v_n})$ is a
pro-$p$-group and $E[l^\infty]^{G_{\F_w}}$ is an $l$-primary group.$\qquad\square$

\begin{rem}\label{cniso}
\emph{ Let $\G_{v_n}:=Gal(\F_w/F_{v_n})$: then the image of $d_w$ is $\G_{v_n}$-invariant. The Hochschild-Serre spectral sequence
\[ H^1(F_{v_n},E[l^\infty])\ri H^1(\F_w,E[l^\infty])^{\G_{v_n}}\ri H^2(\G_{v_n},E[l^\infty](\F_w))=0  \]
shows that $Coker\,d_w=0$ as well. }
\end{rem}

\subsection{The Selmer dual} Theorem \ref{Controll} leads to a (partial) description of the Selmer groups (actually
of their Pontrjagin duals) as modules over the algebra $\Z_l[[\G]]$. Since $\Z_l[[\G]]=\liminv \Z_l[\G/\G_n]$, this
ring is compact with respect to the inverse limit topology. The following generalization of Nakayama's Lemma is
proved in \cite[section 3]{BH}.

\begin{thm}\label{BH}
Let $\L$ be a compact topological ring with 1 and let $J$ be an ideal such that $J^n\rightarrow 0$. Let $X$ be a
profinite $\L$-module. If $X/JX$ is a finitely generated $(\L/J)$-module then $X$ is a finitely generated
$\L$-module.
\end{thm}

\noindent We lack an adequate description of the ideals $J$ of $\Z_l[[\G]]$ such that $J^n\rightarrow 0$. As an
example consider the classical augmentation ideal $I$, i.e. the kernel of the map $\Z_l[[\G]]\ri \Zl$ sending
each $\sigma\in \G$ to 1. It does not verify this condition since $I=I^2$, as the next lemma shows.

\begin{lemma}\label{I=I^2} Let $R$ be a compact topological ring where $p$ is invertible and $G=\liminv G/U_n$,
where the $G/U_n$'s are finite abelian $p$-groups. Denote by $I$ the augmentation ideal of $R[[G]]$: then $I/I^2=0$.
\end{lemma}

\noindent\emph{Proof.} Define $I_n$ as the augmentation ideal of $R[G/U_n]$. The augmentation map on $R[[G]]$
is defined via a limit and $I=\liminv I_n\,$. Besides $I^2=\liminv I_n^2$ because of the compactness hypothesis: the
claim follows if $I_n/I_n^2=0$ for all $n>0$. Let $J_n$ be the augmentation ideal of $\Z[G/U_n]$: since $I_n$ is the free module
generated by the elements $g-1$, we have $I_n=J_n\otimes R$ and $I_n^2=J_n^2\otimes R$. For a commutative group $H$ one has $H\simeq
J_H/J_H^2$ (where $J_H$ is the augmentation of $\Z[H]$): it follows $I_n/I_n^2=(J_n/J_n^2)\otimes
R=0$.$\qquad\square$

Anyway the ideal $lI$ obviously verifies $(lI)^n\rightarrow 0$ and we want to apply Theorem \ref{BH} with
$\L=\Z_l[[\G]]$, $J=lI$ and $X=Sel_E(\F)_l^{\vee}$.

\begin{lemma}\label{lIdual}
Let $M$ be a discrete $\Z_l[[\G]]$-module and $m_l: M\rightarrow M$ the multiplication by $l$. Then
\[ M^{\vee}/lIM^{\vee}\simeq (m_l^{-1}(M^\G))^\vee = (M^\G +M[l])^\vee \]
(where $M[l]$ is the $l$-torsion of $M$).
\end{lemma}

\noindent\emph{Proof.} Let $N$ be the dual of $M$: then $N$ is a $\Z_l[[\G]]$-module. 
Consider the natural projection $\pi : N\sri N/lIN$ and its dual map
$\pi^\vee : (N/lIN)^\vee\iri N^\vee\,$. Let $\phi\in N^\vee\,$: then
\[ \phi\in Im\,\pi^\vee \iff \phi(l(\g -1)\cdot a)=0 \]
for any $\g\in\G$ and any $a\in N$. But $\phi(l\g\cdot a)=\phi(la)$
if and only if $l\phi$ is $\G$-invariant, i.e. $\phi\in Im\,\pi^\vee\iff 
l\phi\in (N^\vee)^\G$. Therefore
\[ (N/lIN)^{\vee}\simeq m_l^{-1}((N^{\vee})^\G)\ .\]
Taking duals one gets
\[  M^{\vee}/lIM^{\vee}\simeq (m_l^{-1}(M^\G))^{\vee}\ .\]
Since $\G$ is pro-$p$ and $M[l]$ is $l$-torsion, one has $H^1(\G,M[l])=0$: therefore
\[ m_l(M)^\G\simeq (M/M[l])^\G\simeq M^\G/M[l]^\G\simeq m_l(M^\G)\ .\]
Thus
\[ m_l^{-1}(M^\G)= m_l^{-1}(m_l(M)^\G)\simeq m_l^{-1}(m_l(M^\G))=M^\G+M[l]\ .\qquad\square \]

\begin{cor}\label{ControllCor}
Assume that both $Sel_E(F)_l$ and $Sel_E(\F)_l[l]$ are finite. Then $Sel_E(\F)_l^{\vee}$ 
is a finitely generated $\Z_l[[\G]]$-module.
\end{cor}

\noindent\emph{Proof.} By the previous lemma with $M=Sel_E(\F)_l$ one has
\[ Sel_E(\F)_l^{\vee}/lI Sel_E(\F)_l^{\vee} \simeq (Sel_E(\F)_l^{\G}+Sel_E(\F)_l[l])^{\vee} \simeq\]
\[ \simeq (Sel_E(F)_l+Sel_E(\F)_l[l])^\vee  \]
(by Theorem \ref{Controll}), so this quotient is finite by hypotheses. Then Theorem \ref{BH} yields
the corollary.$\qquad\square$

In the corollary it would be enough to assume that $Sel_E(F)_l$ and $Sel_E(\F)_l[l]$ are 
cofinitely generated modules over the mysterious ring $\Z_l[[\G]]/lI \Z_l[[\G]]$. Unfortunately 
even with the stronger assumption of finiteness we can't go further (i.e., we are not able to see whether $Sel_E(\F)_l^{\vee}$ is
a torsion $\Z_l[[\G]]$-module or not) due to our lack of understanding of the structure of $\Z_l[[\G]]$-modules
even for simpler $\G$'s like for example $\Zp\,$.

\section{Control theorems for $l=p$} \label{pcontrol}
As stated in the introduction, in this section we shall work with a $\mathbb{Z}_p^d$-extension $\F_d/F$, $d\ge 1$.
As before, it is convenient to write $\F_d$ as the union of finite extensions $F_n$'s of $F$ with $F_n\subset F_{n+1}\,$.\\ 
Moreover we shall consider the $p$-torsion of the elliptic curve $E$:
therefore we use flat cohomology as explained in section \ref{SelGr}, where we described the Selmer groups for this case.
We will mainly follow the notations given there except for the following minor change regarding
the Galois groups: in this section we will write $\G:=Gal(\F_d/F)$ and $\G_n:=Gal(\F_d/F_n)$.

\subsection{Lemmas} We need some lemmas which will be used in the proof of the main theorems.

\begin{lemma}\label{H12Zpd}
Let $\G\simeq \mathbb{Z}_p^d$ and $B$ a finite $p$-primary $\G$-module. Then 
\[ |H^1(\G,B)|\le |B|^d\ {\rm and}\ |H^2(\G,B)|\le |B|^{\frac{d(d-1)}{2}}\ .\]
\end{lemma}

\noindent\emph{Proof.} We use induction on $d$. The case $d=1$ is straightforward since, for $\G=\ov{<\g >}\simeq \Zp\,$,
$H^1(\G,B)\simeq B/(\g-1)B$, so that $|H^1(\G,B)|\le |B|$, and $H^2(\G,B)=0$ because $\Zp$ has $p$-cohomological dimension 1 
(see \cite[3.5.9]{NSW}).\\
For the induction step take $\g$ in a set of independent topological generators of $\G$
and let $\G':=\G/\ov{<\g >}\simeq \Z_p^{d-1}\,$.\\
From the inflation restriction sequence one has 
\[ H^1(\G',B^{\ov{<\g >}})\iri H^1(\G,B)\ri H^1(\ov{<\g >},B)^{\G'} \]
and, using induction,
\[ |H^1(\G,B)|\le |H^1(\G',B^{\ov{<\g >}})||H^1(\ov{<\g >},B)|\le |B|^{d-1}|B|\ .\]
Since $H^i(\ov{<\g >},B)=0$ for all $i\ge 2$, the Hochschild-Serre spectral sequence (see \cite[II.1 Exercise 5]{NSW})
gives an exact sequence
\[ H^2(\G',B^{\ov{<\g >}})\ri H^2(\G,B)\ri H^1(\G',H^1(\ov{<\g >},B))\ . \]
By induction and the result on the $H^1$ one has
\[ |H^2(\G,B)|\le |H^2(\G',B^{\ov{<\g >}})||H^1(\G',H^1(\ov{<\g >},B))|\le \]
\[ \le |B|^{\frac{(d-1)(d-2)}{2}}|B|^{d-1}=|B|^{\frac{d(d-1)}{2}}\ .\qquad\square\]

\begin{lemma}\label{H1(E(L))}
Let $L/K$ be a finite Galois extension of local fields 
and $G$ its Galois group. Let $E/K$ be an elliptic curve
with split multiplicative reduction.
Then $H^1(G,E_0(L))\simeq \Z/e\Z$ (where $e$ is the ramification index
of $L/K$ and $H^1(G,\cdot )$ denotes Galois cohomology).
\end{lemma}

\noindent\emph{Proof.} We recall that Tate parametrization yields an isomorphism of Galois modules 
$E_0(L)\simeq \ol_L^*\,$. The valuation map gives the sequence
\[ \ol_L^*\iri L^*\sri \Z \]
and, via $G$-cohomology,  
\[ \ol^*_K \iri K^*{\buildrel \nu\over\longrightarrow} \Z \sri H^1(G,\ol_L^*) \]
(because $H^1(G,L^*)=0$ by Hilbert 90).
The lemma follows from $Im\,\nu=e\Z$.$\qquad\square$

\begin{lemma}\label{FinpTor}
The group $E[p^\infty](\F_d)$ is finite.
\end{lemma}

\noindent\emph{Proof.} The following argument actually proves that
$E[p^\infty](F^{sep})$ is finite (for a more detailed exposition see \cite[Lemma 2.2]{Br}).\\
Factoring the multiplication-by-$p^m$ map via the $p^m$th power Frobenius one sees that 
$E[p^m](\ov{F})\subset E(F)$ implies $j(E)\in (F^*)^{p^m}\,$.
Therefore if $j(E)\in (F^*)^{p^n}-(F^*)^{p^{n+1}}$ one has 
$E[p^\infty](F^{sep})\subset E[p^n](\ov{F})$.$\qquad\square$

\subsection{The theorem} We are now ready to prove the main theorem: the proof is divided
in several parts and exploits all the techniques which will later lead to similar
results like Theorems \ref{ControlModSel} and \ref{Controlpn}.

\begin{thm}\label{Controlp}
Assume that all primes which are ramified in $\F_d/F$ are of split multiplicative reduction for $E$. 
Then the canonical maps
\[ Sel_E(F_n)_p\longrightarrow Sel_E(\F_d)_p^{\G_n} \]
have finite kernels all bounded by $|E[p^\infty](\F_d)|^d\,$ and cofinitely generated cokernels (over $\Zp\,$). 
\end{thm}

\noindent\emph{Proof.} We start by fixing the notations which will be used throughout the proof.\\ Let
$X_n:=Spec\,F_n\,$, $\X_d:=Spec\,\F_d\,$, $X_{v_n}:=Spec\,F_{v_n}$ and $\X_{w}:=Spec\,\F_{w}\,$ (now $\F_w$ is the
completion of $\F_d$ at $w$).\\
Finally, to ease notations, let
\[ \mathcal{G}(X_n):=Im\big\{H_{fl}^1(X_n,E[p^\infty])\ri
\prod_{v_n\in\mathcal{M}_{F_n}}\,H_{fl}^1(X_{v_n},E[p^\infty])/Im\,\k_{v_n}\big\} \]
(analogous definition for $\mathcal{G}(\X_d)\,$).\\
The map $\X_d\ri X_n$ is a Galois covering with Galois group $\G_n\,$. In this context 
the Hochschild-Serre spectral sequence holds by \cite[III.2.21 a),b) and III.1.17 d)]{Mi1}.
Therefore one has an exact sequence 
\[ \begin{array}{rl} H^1(\G_n,E[p^\infty](\F_d))&\iri H_{fl}^1(X_n,E[p^\infty])\ri  \\ 
\ & \ri H_{fl}^1(\X_d,E[p^\infty])^{\G_n} \ri H^2(\G_n,E[p^\infty](\F_d)) \end{array}\]
which fits into the diagram
\[ \xymatrix{
Sel_E(F_n)_p \ar[d]^{a_n} \ar@{^{(}->}[r] & H_{fl}^1(X_n,E[p^\infty]) \ar[d]^{b_n} \ar@{->>}[r] &
\mathcal{G}(X_n)\ar[d]^{c_n}\\ Sel_E(\F_d)_p^{\G_n} \ar@{^{(}->}[r] & H_{fl}^1(\X_d,E[p^\infty])^{\G_n} \ar[r] &
\mathcal{G}(\X_d) \ .} \] 
As in Theorem \ref{Controll} we shall focus on
\[ Ker\,b_n=H^1(\G_n,E[p^\infty](\F_d))\ ,\
Coker\,b_n\subseteq H^2(\G_n,E[p^\infty](\F_d)) \]
and $Ker\,c_n\,$. 

\subsubsection{$Ker\,b_n$} \label{kerbn} Since $\G_n\simeq\mathbb{Z}_p^d$ and $E[p^\infty](\F_d)$ is finite by
Lemma \ref{FinpTor}, we can apply Lemma \ref{H12Zpd} to get $|Ker\,b_n|\le |E[p^\infty](\F_d)|^d$.

\subsubsection{$Coker\,b_n$} \label{cokerbn} As above we simply use Lemma \ref{H12Zpd} to get 
$|Coker\,b_n|\le |E[p^\infty](\F_d)|^{\frac{d(d-1)}{2}}$.

\begin{rem}\label{kerbnrem}
{\em Note that the result on $Ker\,b_n$ already implies the finiteness of $Ker\,a_n$ and gives
a bound independent of $n$. The bounds for $Ker\,b_n$ and $Coker\,b_n$ depend on $d$ and this is one of the 
reasons we could not consider the full $\Zpn$-extension in this setting. }
\end{rem}

\subsubsection{$Ker\,c_n$}
First of all we note that $Ker\,c_n$ is contained in the kernel of the map
\[  d_n : \prod_{v_n}\,H_{fl}^1(X_{v_n},E[p^\infty])/Im\,\k_{v_n}
\longrightarrow \prod_{w}\, H_{fl}^1(\X_{w},E[p^\infty])/Im\,\k_{w} \] 
and we consider the maps
\[  d_w : H_{fl}^1(X_{v_n},E[p^\infty])/Im\,\k_{v_n}
\longrightarrow H_{fl}^1(\X_{w},E[p^\infty])/Im\,\k_{w} \] 
separately. Note that if $w_1\,,w_2\mid v_n$ then
$Ker\,d_{w_1}\simeq Ker\,d_{w_2}\,$. Moreover, letting 
$d_{v_n}:=\prod_{w|v_n}\,d_w\,$, we have 
\[ Ker\,d_{v_n}= \bigcap_{w|v_n} Ker\,d_w\ {\rm and}\ Ker\,c_n\subseteq
\prod_{v_n\in\mathcal{M}_{F_n}} Ker\,d_{v_n}\ . \]

\subsubsection{Primes of good reduction}\label{GoodRedPr}
Let $v$ be the prime of $F$ lying below $v_n\,$. Assume that $E$ has good reduction
at $v$. In this setting Ulmer proves that
\[ Im(E(F_{v_n})/p^m E(F_{v_n}))\simeq H_{fl}^1(Y_{v_n},E[p^m]) \]
where $Y_{v_n}:=Spec\,\ol_{v_n}$ and $\ol_{v_n}$ is the ring of integers
of $F_{v_n}$ (see \cite[Lemma 1.2]{Ul}).\\
Taking direct limits one finds
\[ Im\,\k_{v_n}= Im(E(F_{v_n})\otimes\Qp/\Zp)\simeq
H_{fl}^1(Y_{v_n},E[p^\infty])\ . \]

\noindent To ease notations let, for any local field $L$,
\[ \mathcal{H}(L):=H_{fl}^1(Spec\,L,E[p^\infty])/H_{fl}^1(Spec\,\ol_L,E[p^\infty])\, . \]
One gets a diagram
\[ \xymatrix {
H_{fl}^1(Y_{v_n},E[p^\infty]) \ar[d]^{\lambda_1} \ar@{^{(}->}[r] & H_{fl}^1(X_{v_n},E[p^\infty]) \ar[d]^{\lambda_2}
\ar@{->>}[r] & \mathcal{H}(F_{v_n}) \ar[d]^{d_w} \\ 
H_{fl}^1(\mathcal{Y}_w,E[p^\infty])^{\G_{v_n}} \ar@{^{(}->}[r] &
H_{fl}^1(\X_w,E[p^\infty])^{\G_{v_n}} \ar[r] & \mathcal{H}(\F_w)^{\G_{v_n}} } \] 
with $\mathcal{Y}_w:=Spec\,\ol_w$ (for the injectivity of the horizontal maps on the left see \cite[III.7]{Mi2}). Note
that the vertical map on the right is exactly $d_w$ because both are induced by the restriction $\lambda_2$ and 
$Im\,d_w$ is $\G_{v_n}$-invariant.\\
The snake lemma yields an exact sequence
\[ Ker\,\lambda_1 \ri Ker\,\lambda_2 \ri Ker\,d_w \ri Coker\,\lambda_1 \ri Coker\,\lambda_2 \ .\]
By our hypothesis all primes of good reduction are unramified: so
for the two maps $\lambda_i$ the situation is identical to the one we described for the map $b_n$, i.e. one has
Galois coverings $\X_w\ri X_{v_n}$ and $\mathcal{Y}_w\ri Y_{v_n}$ both with Galois group $\G_{v_n}:=Gal(\F_w/F_{v_n})$, a subgroup
of $\G_n$ (which one depends on the behaviour of the prime $v_n$ in the extension $\F_d/F_n\,$). Hence
\[ Ker\,\lambda_1=H^1(\G_{v_n},E[p^\infty](\ol_w))\ , \]
\[ Ker\,\lambda_2=H^1(\G_{v_n},E[p^\infty](\F_w))\ , \]
\[ Coker\,\lambda_1\subset H^2(\G_{v_n},E[p^\infty](\ol_w))\ , \]
\[ Coker\,\lambda_2\subset H^2(\G_{v_n},E[p^\infty](\F_w))\ . \]
Since $E$ has good reduction at $v$ (hence at $w$) one has $E(\ol_w)= E(\F_{w})$, therefore $Ker\,\lambda_1\simeq
Ker\,\lambda_2\,$. For the same reason 
\[ H^2(\G_{v_n},E[p^\infty](\ol_w))\simeq H^2(\G_{v_n},E[p^\infty](\F_w)) \]
and the map $Coker\,\lambda_1\ri Coker\,\lambda_2$ (which is induced by this isomorphism) has to be injective.
The snake lemma sequence yields $Ker\,d_w=0$. Note that this is coherent with the number field case, where one
has $Ker\,d_w=0$ for all primes $w$ dividing a prime $l\neq p$ of good reduction (see e.g. \cite[Lemma 4.4]{Gr2}).

\subsubsection{Primes of bad reduction}\label{CTBadRedPr} 
Note that if $v_n$ splits completely then $\F_w= F_{v_n}$ and $d_w$ is clearly an isomorphism.\\
Now let $v_n$ be a prime of bad reduction with nontrivial decomposition group in $\G_n\,$. 
From the Kummer exact sequence we have a diagram
\[ \xymatrix{
H_{fl}^1(X_{v_n},E[p^\infty])/Im\,\k_{v_n}\ar@{^{(}->}[r]\ar[d]^{d_w} & H_{fl}^1(X_{v_n},E)\ar[d]^{h_w} \\
H_{fl}^1(\X_w,E[p^\infty])/Im\,\k_w\ar@{^{(}->}[r] & H_{fl}^1(\X_w,E) } \] 
so that
\[ Ker\,d_w\iri Ker\,h_w\simeq H^1(\G_{v_n},E(\F_w)) \ .\]
Consider the Tate curve exact sequence
\[ q_{E,v}^\Z\iri \F_w^*\sri E(\F_w) \]
and take Galois cohomology to get
\[ H^1(\G_{v_n},E(\F_w))\iri H^2(\G_{v_n},q_{E,v}^\Z) \]
where the injectivity comes from Hilbert 90.\\
Since $q_{E,v}\in F_v$ the action of $\G_{v_n}$ on $q_{E,v}^\Z$ is trivial, hence
\[ Ker\,d_w\iri H^2(\G_{v_n},q_{E,v}^\Z)\simeq H^2(\G_{v_n},\Z)\simeq (\G_{v_n}^{ab})^\vee\simeq(\Qp/\Zp)^{d(v_n)} \]
where $d(v_n)={\rm rank}_\Zp\,\G_{v_n}\le d$ (see \cite[pag.\,50]{NSW}).\\
Therefore $Ker\,d_w$ and $Ker\,d_{v_n}=\bigcap Ker\,d_w$ are cofinitely generated and,
since there are only finitely many $v_n$'s of bad reduction, $Ker\,c_n$ is cofinitely
generated as well.$\qquad\square$

\begin{rem}\label{ControlpRem}
\emph{In the main examples we have in mind for $\F/F$ (see sections \ref{pLfunc} and \ref{pLf}) there is only one
prime which ramifies and it totally ramifies in $\F/F$. Moreover the available constructions for $p$-adic
$L$-functions require this prime to be of split multiplicative reduction so our Theorem \ref{Controlp}
applies to these settings.}
\end{rem}

\begin{rem}\label{RankRem}
{\em Note that, due to the (possibly) increasing number of primes of bad reduction in the tower of
$F_n$'s, the coranks of the $Ker\,c_n$'s are finite but not bounded in general. For example if $v$ is 
inert unramified in $\F_d/F$ and $d\ge 2$ then there are infinitely many primes in $\F_d$ lying over $v$.}
\end{rem}

\noindent Recall that $\L_d$ acts on $Sel_E(\F_d)_p\,$ and on $\S_d:=Sel_E(\F_d)_p^{\vee}\,$.

\begin{cor}\label{ControlpCor} In the setting of Theorem \ref{Controlp}, assume that $Sel_E(F)_p$ is a
cofinitely generated $\Zp$-module. Then $\S_d$ is a finitely generated $\L_d$-module.
\end{cor}

\noindent\emph{Proof.} Let $I$ be the augmentation ideal of $\L_d\,$.
One has that $\S_d/I\S_d$ is dual to $Sel_E(\F_d)_p^\G$
which, by the hypothesis on $Sel_E(F)_p$ and Theorem \ref{Controlp}, is cofinitely generated over $\Zp$.
Therefore $\S_d$ is a finitely generated $\L_d$-module by Theorem \ref{BH}$\,$.$\qquad\square$

\subsection{Primes of bad reduction}\label{BadRedPr}
We can be more precise on $Ker\,d_w$ for the inert primes of bad reduction.

\begin{lemma}\label{BadRedPrLemma}
With the hypothesis and notations of Theorem \ref{Controlp}, let $v_n$ be a place of
bad reduction. One has:
\begin{itemize}
\item[{\bf 1.}] $Ker\,d_w$ is finite if $v_n$ is unramified in $\F_d/F$;
\item[{\bf 2.}] $Ker\,d_w=0$ if $v_n$ splits completely in $\F_d/F$;
\item[{\bf 3.}] $corank_\Zp Ker\,d_w\le d$ if $v_n$ is ramified in $\F_d/F$.
\end{itemize}
\end{lemma}

\noindent\emph{Proof.} Part {\bf 2} and {\bf 3} have been proven in section \ref{CTBadRedPr}.
For the other primes consider the embedding coming from the Kummer exact sequence
\[ Ker\,d_w\iri Ker\,h_w\simeq H^1(\G_{v_n},E(\F_w))\simeq
\dl{m} H^1(\G_{v_n}^{v_m},E(F_{v_m})) \]
where $\G_{v_n}^{v_m}:=Gal(F_{v_m}/F_{v_n})$ is a quotient of $\G_{v_n}\subset \Z_p^d\,$.
Consider the exact sequence
\[ E_0(F_{v_m})\iri  E(F_{v_m})\sri T_{v_m} \]
where $T_{v_m}$ is a cyclic group of order $-ord_{v_m}(j(E))$. Our Lemma \ref{H1(E(L))} applies here with $L=F_{v_m}$
and $K=F_{v_n}$ and one gets 
\[ T_{v_m}^{\G_{v_n}^{v_m}}/T_{v_n}\iri \Z/e_{v_m/v_n}\Z \rightarrow
H^1(\G^{v_m}_{v_n},E(F_{v_m}))\rightarrow H^1(\G^{v_m}_{v_n},T_{v_m}) \]
where $e_{v_m/v_n}$ is the ramification index of $F_{v_m}/F_{v_n}\,$, hence a power of $p$.
From Tate parametrization one has an isomorphism of Galois modules
$ T_{v_m}\simeq F_{v_m}^*/\ol_{v_m}^*q_{E,v_m}^\Z\,$. 
Let $\pi_{v_m}$ be a uniformizer for the prime $v_m\,$:
then for any $\sigma\in \G^{v_m}_{v_n}$ one has $\sigma(\pi_{v_m})=
u\pi_{v_m}$ for some $u\in \ol_{v_m}^*\,$. Hence $\G^{v_m}_{v_n}$ acts trivially
on $T_{v_m}$ and
\[ H^1(\G^{v_m}_{v_n},T_{v_m})\iri (T_{v_m,p})^{d(v_n)} \]
(where $T_{v_m,p}$ is the $p$-part of the group $T_{v_m}$ and $d(v_n)$ is the number of generators of $\G_{v_n}^{v_m}\,$). 
If $v_n$ is unramified then, for any $m$, $e_{v_m/v_n}=1$, $d(v_n)\le 1$ and 
\[ |T_{v_m}|=| T_{v_n}|=|T_v|=-ord_v(j(E)) \] 
is finite and constant, so $H^1(\G_{v_n},E(\F_w))\iri T_{v,p}\,$. 
Thus $Ker\,d_w$ is finite and this proves {\bf 1}$\,$.$\qquad\square$

\begin{rem}\label{BadRedPrRemRam}
\emph{With the notations of the previous lemma, if $v_n$ is ramified then
\[ |T_{v_m}|=-ord_{v_m}(j(E))=-ord_{v_n}(j(E))e_{v_m/v_n}=|T_{v_n}|e_{v_m/v_n} \ .\]
Since the $T_{v_k}$'s are cyclic and the action of $\G_{v_n}^{v_m}$ is trivial one gets 
\[ T_{v_m}^{\G_{v_n}^{v_m}}/T_{v_n} \simeq \Z/e_{v_m/v_n}\Z \]
and one finds injections
\[ H^1(\G_{v_n}^{v_m},E(F_{v_m}))\iri H^1(\G_{v_n}^{v_m},T_{v_m})\iri (T_{v_m,p})^{d(v_n)} \ .\]
Hence, taking limits, one has
\[ H^1(\G_{v_n},E(\F_w))\iri \dl{m} (T_{v_m,p})^{d(v_n)}\simeq (\Qp/\Zp)^{d(v_n)} \]
which gives another proof of part {\bf 3} of Lemma \ref{BadRedPrLemma}.}
\end{rem}

\subsection{Modified Selmer groups}
In section \ref{algL} we are going to consider characteristic and Fitting ideals of
$Sel_E(\F_d)_p^\vee$ in order to define some algebraic $L$-function for $E$. Those ideals
are well defined for the finitely generated modules provided by Corollary \ref{ControlpCor},
but they are non-trivial only for torsion modules. In order to reduce to this latter case one needs 
finite kernels and cokernels for the maps
\[ Sel_E(F_n)_p\longrightarrow Sel_E(\F_d)_p^{\G_n} \]
and in Theorem \ref{Controlp} and Lemma \ref{BadRedPrLemma} we have seen that the only 
obstruction comes from the ramified primes. We are going to modify the definition of the 
Selmer groups in order to give at least two ways in which one could get torsion modules.

\begin{defin}\label{DefModSel}
With all notations as in section \ref{Selmer} let $\Sigma\subset \mathcal{M}_L$ be a set of places of an algebraic 
extension $L$ of $F$. We define 
\[ Sel_{E,\Sigma}(L)_p :=Ker\{ H_{fl}^1(X_L,E[p^\infty])\ri\prod_{v\not\in\Sigma}
H_{fl}^1(X_{L_v},E[p^\infty])/Im\,\k_v \,\} \]
and
\[ Sel_E^\Sigma(L)_p := Ker\{ H_{fl}^1(X_L,E[p^\infty]) - \]
\[  \ri \prod_{v\not\in\Sigma} 
H_{fl}^1(X_{L_v},E[p^\infty])/Im\,\k_v \,\times\,\prod_{v\in\Sigma}H_{fl}^1(X_{L_v},E[p^\infty])\,\} \ . \]
to be the $p$-parts of the (lower and upper respectively) \emph{$\Sigma$-Selmer groups} of $E$ over $L$.
\end{defin}

Note that  
\[ Sel_E^\Sigma(L)_p \subseteq Sel_E(L)_p\subseteq Sel_{E,\Sigma}(L)_p \]
with equality occurring if $\Sigma=\emptyset$.\\
In the same setting of Theorem \ref{Controlp} we can now prove

\begin{thm}\label{ControlModSel}
Let $\F_d/F$ be a $\Z_p^d$-extension. Let $\Sigma(d)$ be the set of ramified places in $\F_d$ and, 
for any $n$, let $\Sigma_n$ be the set of primes lying below $\Sigma(d)$. 
Assume that $E$ has split multiplicative reduction at all places in $\Sigma_0$ (so that, in particular,
$\Sigma_0$ is finite). Then the maps
\[ \alpha_n\,:\, Sel_{E,\Sigma_n}(F_n)_p\rightarrow Sel_{E,\Sigma(d)}(\F_d)_p^{\G_n} \]
and
\[ \beta_n\,:\, Sel_E^{\Sigma_n}(F_n)_p\rightarrow Sel_E^{\Sigma(d)}(\F_d)_p^{\G_n} \]
have finite kernels and cokernels.
\end{thm}

\noindent\emph{Proof.} The structure of the proof is the same as for Theorem \ref{Controlp}: no changes are
needed for the maps $b_n$ and $d_w$'s for all $w$'s of good reduction. We are left with the finitely many
places of bad reduction in $F_n\,$. For unramified primes 
Lemma \ref{BadRedPrLemma} shows that $Ker\,d_w$ is finite and this completes the proof for the maps $\alpha_n\,$.\\
For the upper $\Sigma$-Selmer groups one has to consider also the maps
\[ \t d_w\,:\,H^1_{fl}(X_{v_n},E[p^\infty])\rightarrow H^1_{fl}(X_w,E[p^\infty]) \]
for the ramified primes $v_n\,$. As seen before one can apply Hochschild-Serre to find
$Ker\,\t d_w=H^1(\G_{v_n},E[p^\infty](\F_w))$ and, since $E[p^\infty](\F_w)$ is finite (same proof as
Lemma \ref{FinpTor}, using $\F_w\bigcap \ov{F} \subseteq F^{sep}\,$), one gets
\[ |Ker\,\t d_w|\le |E[p^\infty](\F_w)|^d \]
by Lemma \ref{H12Zpd}.$\qquad\square$

\begin{cor}\label{ControlModSelCor}
In the setting of Theorem \ref{ControlModSel} assume that $Sel_{E,\Sigma_0}(F)_p$ is cofinitely generated 
over $\Zp$ (resp. finite). Then $Sel_{E,\Sigma(d)}(\F_d)_p$ is a cofinitely generated (resp. torsion) 
$\L_d$-module. The same statement holds (with identical hypotheses on $Sel_E^{\Sigma_0}(F)_p\,$)
for the upper $\Sigma(d)$-Selmer group. 
\end{cor}

\noindent\emph{Proof.} The proof of Corollary \ref{ControlpCor} applies here as well.
One only has to note that if $M/IM$ is finite for some finitely generated $\L_d$-module $M$
(where $I$ is the augmentation ideal of $\L_d$) then $M$ is a torsion module by the final Theorem
of \cite[section 4]{BH}.$\qquad\square$

For a $\Zp$-extension we can therefore give a proof of the following corollary
along the lines of \cite[Corollary 4.9 and Theorem 1.3]{Gr2}.

\begin{cor}\label{ControlModSelCor2}
In the setting of the previous corollary assume that $Sel_{E,\Sigma_0}(F)_p$ is finite. 
Then $E(\F_1)$ is finitely generated.
\end{cor}

\noindent\emph{Proof.} Let $\Sigma:=\Sigma(1)$ and let $\S_{1,\Sigma}$ be the dual of $Sel_{E,\Sigma}(\F_1)_p\,$.
By Corollary \ref{ControlModSelCor} $\S_{1,\Sigma}$ is a finitely generated torsion $\L_1$-module. 
By the well-known structure theorem for such modules there is a pseudo-isomorphism
(i.e. with finite kernel and cokernel)
\[ \S_{1,\Sigma}\sim \bigoplus_{i=1}^r \Zp[[T_1]]/(f_i^{e_i})\,. \]
Let $\lambda=\deg \prod f_i^{e_i}\,$: then ${\rm rank}_\Zp \S_{1,\Sigma}=\lambda$ and, taking duals, one gets
\[ (Sel_{E,\Sigma}(\F_1)_p)_{div}\simeq (\Qp/\Zp)^\lambda \]
(where $(Sel_{E,\Sigma}(\F_1)_p)_{div}$ is the divisible part of $Sel_{E,\Sigma}(\F_1)_p\,$). 
By Theorem \ref{ControlModSel}, for any $n$, one has
\[ (Sel_{E,\Sigma_n}(F_n)_p)_{div}\simeq (\Qp/\Zp)^{t_n}\ {\rm with}\
t_n\le\lambda\ . \]
Moreover we know that
\[ E(F_n)\otimes \Qp/\Zp \simeq (\Qp/\Zp)^{rank\,E(F_n)} \]
and, obviously,
\[ E(F_n)\otimes \Qp/\Zp \iri (Sel_{E,\Sigma_n}(F_n)_p)_{div} \]
(which is not true in general for the upper $\Sigma$-Selmer groups).
Therefore $rank\,E(F_n)\le t_n\le \lambda$ for any $n$, i.e. such ranks are bounded.\\ 
Choose $m$ such that $rank\,E(F_m)$ is maximal, then $E(\F_1)/E(F_m)$ is a torsion group. 
Take $P\in E(\F_1)$ and let $s$ be such that $sP\in E(F_m)$. Then for any $\g\in Gal(\F_1/F_m)$ 
one has $s(\g(P)-P)=0$, i.e. $\g(P)-P\in E(\F_1)_{tor}$. Since the torsion points in $E(\F_1)$ 
are finite (we provide a proof in Lemma \ref{torsion} below) take
$t=|E(\F_1)_{tor}|$ to get $t(\g(P)-P)=0$. Thus $tP\in
E(F_m)$ for all $P\in E(\F_1)$ and multiplication by $t$ gives a homomorphism $\varphi_t:\ E(\F_1)\rightarrow
E(F_m)$ whose image is finitely generated (being a subgroup of $E(F_m)\,$) and whose kernel is the finite group
$E(\F_1)_{tor}\,$. Hence $E(\F_1)$ is indeed finitely generated.$\qquad\square$

\begin{lemma} \label{torsion} 
For any $d\geq 1$, the set $E(\F_d)_{tor}$ is finite. 
\end{lemma}

\noindent\emph{Proof.} For any $l$, let $K_l$ be the minimal extension of $F$ such that 
$E[l](K_l)=E[l](\ov F)$. By Igusa's work \cite{Ig}, it is known that for almost all primes 
$l\neq p$ the Galois group $Gal(K_l/F)$ contains a subgroup isomorphic to $SL_2(\mathbb F_l)$. 
This implies that the subgroup generated by the Galois orbit of any $P\in E[l](\ov F)-\{0\}$
is all of $E[l](\ov F)$. In particular, since $\F_d/F$ is Galois, $E[l](\F_d)\neq\{0\}$ 
would imply $K_l\subset\F_d$, contradicting the fact that $\F_d/F$ is abelian.\\
We are left with the possibility that $E[l^{\infty}](\F_d)$ is infinite for some prime $l$. 
Assume that this happens: then one can choose an infinite 
sequence $P_n\in E[l^n]$ so that $lP_{n+1}=P_n$. Let $\F'\subset\F_d$ be the minimal 
extension of $F$ such that $\{P_n\}_{n\in\N}\subset E(\F')$ and put $\G':=Gal(\F'/F)$. Then $\G'$ is 
both an infinite subgroup of $\Z_l^{\ast}$ (the automorphisms of the group 
generated by the $P_n$'s) and a quotient of $\G\simeq\Z_p^d$: this is impossible for $l\neq p$. We 
already proved that $E[p^{\infty}](F^{sep})$ is finite in Lemma \ref{FinpTor}.$\qquad\square$

\section{The algebraic $L$-function and the main conjecture} \label{algL}

\subsection{The characteristic and Fitting ideals}\label{charfitt} 
In this section, we work again with a $\Zpn$-extension 
$\F/F$. In order to apply results from section \ref{pcontrol}, we consider $\Z_p^d$-subextensions
$\F_d\subset \F$, naturally ordered by inclusion; we call such fields \emph{$\Zp$-finite extensions} (of $F$). 
We will let $\F_d$ vary among all $\Z_p$-finite subextensions of $\F$: 
therefore we need to refine our notations.\\
Let $\L(\F_d):= \Zp[[Gal(\F_d/F)]]$ and $\S(\F_d):=Sel_E(\F_d)_p^{\vee}$ 
(shortened into $\L_d$ and $\S_d$ when it is clear to which extension they refer).

It is assumed throughout that all places which ramify in $\F/F$ are of split multiplicative reduction
for $E$ and that $Sel_E(F)_p$ is cofinitely generated as a $\Zp$-module.\\ 
For any fixed $\F_d\,$, Corollary \ref{ControlpCor} states that $\S_d$ is a 
finitely generated $\L_d$-module. Recall that $\L_d$ is (noncanonically) isomorphic to $\Zp[[T_1,\dots ,T_d\,]]$.
A finitely generated torsion $\L_d$-module is said to be
{\em pseudo-null} if it has at least two relatively prime annihilators (or, equivalently, if its
annihilator ideal has height at least 2). For example if $d=1$ pseudo-null is
equivalent to finite. A {\em pseudo-isomorphism} between finitely generated $\L_d$-modules 
$M$ and $N$ (i.e. a morphism with pseudo-null kernel and cokernel) will be denoted by $M\sim_{\L_d}N$. 
If $M$ is a finitely generated $\L_d$-module then there is a pseudo-isomorphism
\[ M\sim_{\L_d} \L_d^r \oplus \bigg(\bigoplus_{i=1}^{n(M)} \L_d/(g_i^{e_i})\L_d \bigg) \]
where the $g_i$'s are irreducible elements of $\L_d$ (determined up to an element of $\L_d^*\,$) 
and $r$, $n(M)$ and the $e_i$'s are uniquely determined (see e.g. the structure theorem
\cite[VII.4.4 Theorem 5]{Bo}).

\begin{defin}\label{defcharid}
In the above setting the {\em characteristic ideal} of $M$ is
\[ Char_{\L_d}(M):=\left\{ \begin{array}{ll}
0 & {\rm if\ }M\ {\rm is\ not\ torsion} \\
\left(\prod_{i=1}^{n(M)} g_i^{e_i}\right) & {\rm otherwise}
\end{array} \right. \ .\]
\end{defin}

Let $Z$ be a finitely generated $\L_d$-module and let 
\[ \L_d^a {\buildrel\varphi\over\longrightarrow} \L_d^b \sri Z \]
be a presentation where the map $\varphi$ can be represented by an $a\times b$ matrix $\Phi$
with entries in $\L_d\,$.

\begin{defin}\label{defFittid}
In the above setting the {\em Fitting ideal} of $Z$ is
\[ Fitt_{\L_d}(Z):=\left\{ \begin{array}{ll}
0 & {\rm if}\ a<b \\
{\rm the\ ideal\ generated\ by\ all\ the}& \\
{\rm determinants\ of\ the\ } b\times b & {\rm if}\ a\ge b \\
{\rm minors\ of\ the\ matrix\ } \Phi & 
\end{array} \right. \ .\] 
\end{defin}

For the basic theory and properties of Fitting ideals the reader is referred to the Appendix in 
\cite{MW}. Here we only mention the fact that $Fitt_{\L_d}(Z)$ is independent from the presentation
and that if $Z$ is an elementary module, i.e. if
\[ Z = \L_d^r \oplus \bigg(\bigoplus_{i=1}^s \L_d/(g_i^{e_i})\L_d \bigg) \]
then $Fitt_{\L_d}(Z)=Char_{\L_d}(Z)$.

\subsubsection{$\Zp$-finite extensions} We can extend our control theorem to a relation between Selmer groups of
$\Zp$-finite extensions of $F$. 

\begin{thm}\label{Controlpn}
For any inclusion of $\Zp$-finite extensions $\F_d\subset \F_e\,$, $e > d\ge 2$, one has 
\[ Sel_E(\F_d)_p\sim_{\L_d}Sel_E(\F_e)_p^{\G^e_d} \]
(where $\G^e_d=Gal(\F_e/\F_d)\,$).
\end{thm}

\noindent{\em Proof.} As in the proof of Theorem \ref{Controlp}, consider the diagram
\[ \xymatrix{
Sel_E(\F_d)_p \ar[d]^a \ar@{^{(}->}[r] & H_{fl}^1(\X_d,E[p^\infty]) \ar[d]^b \ar@{->>}[r] &
\mathcal{G}(\X_d)\ar[d]^c\\ 
Sel_E(\F_e)_p^{\G^e_d} \ar@{^{(}->}[r] & H_{fl}^1(\X_e,E[p^\infty])^{\G^e_d} \ar[r] &
\mathcal{G}(\X_e) \ .} \]
By Corollary \ref{ControlpCor} we know that $Ker\,a$ and $Coker\,a$ are cofinitely generated $\L_d$-modules.
Moreover $Ker\,a$, $Ker\,b$ and $ Coker\,b$ are finite. 
Therefore it is enough to find two relatively prime annihilators for $Ker\,c$.
We consider the maps
\[ d_{w_e}:H_{fl}^1(\X_{w_d},E[p^\infty])/Im\,\k_{w_d}\ri H_{fl}^1(\X_{w_e},E[p^\infty])/Im\,\k_{w_e} \]
(for any $w_e\in \mathcal{M}_{\F_e}$ dividing $w_d\in\mathcal{M}_{\F_d}\,$).
As in section \ref{GoodRedPr} one sees that $Ker\,d_{w_e}=0$ for primes of good reduction and for
primes which split completely in $\F_e/\F_d\,$.\\
For primes $w_d$ of bad reduction which do not split completely, working as in Theorem \ref{Controlp},
one finds injections
\[ Ker\,d_{w_e}\iri H^1(\G^{w_e}_{w_d},E(\F_{w_e}))\iri H^2(\G^{w_e}_{w_d},q_{E,v}^\Z) \]
(where $\G^{w_e}_{w_d}$ is the local Galois group).\\
The group $Gal(\F_d/F)$ acts trivially both on $q_{E,v}\in F_v$ and on $\G^{w_e}_{w_d}$ (because $\G$ is abelian).
Since $d\ge 2$ there are two topologically independent elements $\g_1\,,\g_2$ of $Gal(\F_d/F)$.
Hence 
\[ H^2(\G^{w_e}_{w_d},q_{E,v}^\Z)\sim_{\L_d}0 \]
because its annihilator ideal contains at least $\g_1-1$ and $\g_2-1$. Therefore all $Ker\,d_{w_e}$'s
are pseudo-null $\L_d$-modules and have common annihilators which then annihilate $Ker\,c$
as well.$\qquad\square$

\begin{cor}\label{ControlpnCor}
In the same setting of Theorem \ref{Controlpn} let $\pi^e_d:\L_e\ri \L_d$ be the canonical projection
and $I_{e/d}$ its kernel. Then
\[ Char_{\L_d}(\S_d)=Char_{\L_d}(\S_e/I_{e/d}\S_e) \ .\]
\end{cor}

\noindent\emph{Proof.} Just take duals in the pseudo-isomorphism given 
by Theorem \ref{Controlpn}.$\qquad\square$

\begin{lemma}\label{PrLimFitt}
Let $\F_d\subset \F_e$ be an inclusion of $\Zp$-finite extensions, $e > d\ge 2$. 
Assume that $E[p^\infty](\F)=0$ or that $Fitt_{\L_d}(\S_d)$ is principal. Then 
\[ \pi^e_d(Fitt_{\L_e}(\S_e))\subseteq Fitt_{\L_d}(\S_d) \ .\]
\end{lemma}

\noindent\emph{Proof.} Dualising in Theorem \ref{Controlpn} one finds a sequence
\[ \S_e/I_{e/d}\S_e \ri \S_d \sri (Ker\,a)^\vee \]
where 
\[ Ker\,a \iri H^1(\G_d^e,E[p^\infty](\F_e)) \]
is finite (by Lemmas \ref{FinpTor} and \ref{H12Zpd}). If $E[p^\infty](\F)=0$ then 
$(Ker\,a)^\vee=0$ and
\[ \pi^e_d(Fitt_{\L_e}(\S_e))= Fitt_{\L_d}(\S_e/I_{e/d}\S_e)\subseteq 
Fitt_{\L_d}(\S_d) \]
by \cite[properties 1 and 4]{MW}.\\
If $E[p^\infty](\F)\neq 0$, by \cite[property 9]{MW} one gets
\[ Fitt_{\L_d}(\S_e/I_{e/d}\S_e)Fitt_{\L_d}((Ker\,a)^\vee) \subseteq 
Fitt_{\L_d}(\S_d) \ .\]
The Fitting ideal of a finitely generated torsion module contains a power of its annihilator
(\cite[property 8]{MW}) so let $\sigma_1\,,\sigma_2$ be two relatively prime elements of
$Fitt_{\L_d}((Ker\,a)^\vee)$. (Recall that $\L_d$ is a unique factorization domain and relatively
prime means that they have no common factor.) Let $\theta_d$ be a generator of $Fitt_{\L_d}(\S_d)$. Then
for any $\alpha\in Fitt_{\L_d}(\S_e/I_{e/d}\S_e)$ one finds that $\sigma_1\alpha$ and
$\sigma_2\alpha$ are divisible by $\theta_d$ (or $\sigma_1\alpha=\sigma_2\alpha=0$ if $\theta_d=0$). 
Hence $\theta_d$ divides $\alpha$ (or $\alpha=0$) and
\[  \pi^e_d(Fitt_{\L_e}(\S_e))=Fitt_{\L_d}(\S_e/I_{e/d}\S_e)
\subseteq Fitt_{\L_d}(\S_d) \ .\qquad\square \]

\begin{rem}\label{PrLimFittRem}
\emph{The hypothesis on $Fitt_{\L_d}(\S_d)$ is verified, for example, if $\S_d$ is an elementary
module or if one can find a presentation with the same number of generators and relations.
Also $E[p^\infty](\F)=0$ if $j(E)\not\in (F^*)^p$ (see Lemma \ref{FinpTor}).}
\end{rem}

\subsubsection{Going to the limit}\label{LimFitt}
For $\F_d\subset\F$ a $\Zp$-finite extension of $F$ we let $\pi_{\F_d}:\L\ri \L(\F_d)$ be the canonical projection.
All $\L(\F_d)$-modules can be thought of as modules over the ring $\L$ with trivial action of
$Gal(\F/\F_d)$.\\
Define
\[ \wt{Fitt}_\L(\S(\F_d)):= \pi_{\F_d}^{-1} (Fitt_{\L(\F_d)}(\S(\F_d))) \ .\]
 
\begin{lemma}\label{}
Under the same assumptions of Lemma \ref{PrLimFitt} one has
\[ \wt{Fitt}_\L(\S_e) \subseteq \wt{Fitt}_\L(\S_d) \ .\]
\end{lemma}

\noindent\emph{Proof.} By Lemma \ref{PrLimFitt} 
\[ \pi^e_d(Fitt_{\L_e}(\S_e))\subseteq Fitt_{\L_d}(\S_d) \ .\]
The claim follows observing that $\pi_{\F_d}=\pi^e_d\circ\pi_{\F_e}$.$\qquad\square$

\begin{defin}\label{DefFittId} With the assumptions of Lemma \ref{PrLimFitt}, 
let $\S:=Sel_E(\F)_p^{\vee}\,$. Its \emph{pro-Fitting ideal} is  
\[ \wt{Fitt}_\L(\S):=\bigcap \wt{Fitt}_\L(\S(\F_d)) \]
where the intersection is taken over all $\Zp$-finite subextensions.
\end{defin}

\begin{rem}\label{DefFittIdRem}
1. \emph{By Corollary \ref{ControlpCor} the $\L$-modules $\S(\F_d)$ are finitely generated: hence one can define their Fitting ideals 
$Fitt_\L(\S(\F_d))$ in the usual way, as the ideals generated by maximal minors of the matrix expressing relations between a (finite) chosen set of generators. 
Then, by \cite[property 4]{MW},
\[ \wt{Fitt}_\L(\S(\F_d))=Fitt_\L(\S(\F_d))+Ker(\pi_{\F_d}) \ .\] 
When also $\S$ is finitely generated one constructs $Fitt_\L(\S)$ in the same way.
Moreover if $E[p^\infty](\F)=0$, then the natural maps $\S\ri \S(\F_d)$ are surjective and 
\[  Fitt_\L(\S)=\wt{Fitt}_\L(\S) \]
by \cite[property 10]{MW} (the result is true also without the noetherian hypothesis in the reference).}\\
2. \emph{ Since $\L=\liminv \L(\F_d)$, one has an equivalent definition using projective limits. Indeed
\[ \wt{Fitt}_\L(\S)=\bigcap_{\F_d} \pi_{\F_d}^{-1} (Fitt_{\L(\F_d)}(\S(\F_d))) = \]
\[ = \il{\F_d} Fitt_{\L(\F_d)}(\S(\F_d)) \ .\]
In particular if $\S$ is finitely generated and $E[p^\infty](\F)=0$, then
\[ Fitt_\L(\S)=\wt{Fitt}_\L(\S) = \il{\F_d} Fitt_{\L(\F_d)}(\S(\F_d)) \]
(for a similar result on the behaviour of Fitting ideals with respect to projective limits
see \cite[Theorem 2.1]{GK}). }
\end{rem}

The use of Fitting ideals instead of characteristic ideals is justified exactly by their behaviour
with respect to limits and some reformulations of the Main Conjectures of Iwasawa theory in terms
of Fitting ideals have already been given in \cite{Gr}, \cite{GK} and \cite{Ku}. Anyway the two ideals 
are strictly connected as the following lemma shows. The ideal $\wt{Fitt}_\L(\S)$ 
(or, if it is principal, a generator of such an ideal) might be a good candidate for our algebraic
$L$-function.

\begin{lemma}\label{FittChar}
Let $\F_d$ be a $\Z_p^d$-extension and $\L_d:=\L(\F_d)$.
Let $M$ be a finitely generated $\L_d$-module. Then 
\[ Fitt_{\L_d}(M)\subseteq Char_{\L_d}(M) \ .\] 
Moreover if $Fitt_{\L_d}(M)$ is principal then $Fitt_{\L_d}(M) = Char_{\L_d}(M)$.
\end{lemma}

\noindent\emph{Proof.} If $M$ is not torsion then $Fitt_{\L_d}(M) = Char_{\L_d}(M)=0$
so we can assume that $M$ is a torsion $\L_d$-module.\\
Let $E$ be the elementary module which is pseudo-isomorphic to $M$ by the structure theorem stated before
Definition \ref{defcharid}. Being pseudo-isomorphic is an equivalence relation between finitely generated torsion
$\L_d$-modules so the structure theorem gives rise to two exact sequences
\[ A_1\iri M \ri E \sri A_2 \]
and
\[ B_1\iri E \ri M \sri B_2 \]
where $A_1\,,A_2\,,B_1\,,B_2$ are pseudo-null. From these sequences one gets
\[ Fitt_{\L_d}(M)Fitt_{\L_d}(A_2)\subseteq Fitt_{\L_d}(E)=Char_{\L_d}(E)=Char_{\L_d}(M) \]
(by \cite[property 9]{MW}, using $Fitt_{\L_d}(M)\subseteq Fitt_{\L_d}(M/A_1)\,$) and
\[ Char_{\L_d}(M)Fitt_{\L_d}(B_2)=Fitt_{\L_d}(E)Fitt_{\L_d}(B_2)\subseteq Fitt_{\L_d}(M)\ . \]
As seen in Lemma \ref{PrLimFitt} let $\sigma_1\,,\sigma_2$ be two relatively prime elements
of $Fitt_{\L_d}(B_2)$ and $\tau_1\,,\tau_2$ two relatively prime elements of $Fitt_{\L_d}(A_2)$.
Let $\theta$ be a generator of $Char_{\L_d}(M)$: then, for any $\alpha\in Fitt_{\L_d}(M)$,
$\alpha\tau_1$ and $\alpha\tau_2$ are divisible by $\theta\,$. Hence $\theta$ divides $\alpha$ and
$Fitt_{\L_d}(M)\subseteq Char_{\L_d}(M)$.\\  
If $Fitt_{\L_d}(M)$ is principal we can use the same proof with the $\sigma_i$'s in the place of the $\tau_i$'s
to get the reverse inclusion and eventually the equality $Fitt_{\L_d}(M) = Char_{\L_d}(M)$.$\qquad\square$

\subsection{The analytic side}\label{pLf} We briefly describe how to associate $p$-adic 
$L$-functions (in the sense of \ref{pLfunc}) to $E$ and to certain Galois extensions $\wt\F/F$.
Since our goal is just to provide an introduction to a main conjecture, the angel of brevity 
compels us to be very sketchy; for the missing details the reader is referred to the original 
papers and to \cite{BL}. 

All examples we know can be seen as applications of the following general ideas.\\
To begin with, we fix a place $\infty$ such that $E$ has conductor $\mathfrak n\infty$ (in particular, because of our 
initial assumption, $E$ has split multiplicative reduction at $\infty$). Then, thanks to the analytic uniformization 
of elliptic curves by Drinfeld modular curves, one can associate to $E$ a $\Z$-valued measure 
$\mu_{E,\infty}$ on $\mathbb{P}^1(F_\infty)$ (\cite[page 386]{Lo}). The next ingredients are a quadratic algebra $K/F$ (the case
$K=F\times F$ is allowed) and an embedding $\Psi:K\ri M_2(F)$. The algebra $K$ and the map 
$\Psi$ are required to satisfy certain conditions which it would be too long to discuss here.\\ 
At this point there are essentially two constructions. One (carried out in \cite{Lo}) consists 
in taking a fundamental domain $X$ for the action on $\mathbb{P}^1(F_\infty)$ of a 
certain quotient of the group $(K\otimes F_{\infty})^*$ via $\Psi$.\\
In the second approach, developed in \cite{HL}, one chooses a second place $\p$ dividing 
$\mathfrak n$. Then, by a construction reminiscent of modular symbols, the fixed points $x,y$ of $\Psi(K^*)$
and the measure $\mu_{E,\infty}$ (or rather the harmonic cocycle attached to it) are employed
to define a measure $\mu_E\{x\ri y\}$ on $\mathbb{P}^1(F_{\p})$.\\
In both cases, class field theory allows to identify the set $X$ (respectively 
$\mathcal O_{\p}^*\subset\mathbb{P}^1(F_{\p})\,$)
with a Galois group $G:=Gal(\wt\F/H)$, where $\wt\F$ is either an anticyclotomic extension of 
$K$ (if $K$ is a field) or a cyclotomic extension of $F$ and $H/F$ is a finite extension. We 
skip precise definitions, just remarking that $G$ contains a finite index subgroup isomorphic to $\Zpn$.\\ 
The restriction of $\mu_{E,\infty}$ to $X$ (resp. of $\mu_E\{x\ri y\}$ to $\mathcal O_{\p}^*$) can be 
thought of as a measure on $G$, that is, an element $L(E)\in \Z[[G]]\subset \Zp[[G]]$. Teitelbaum's measure 
referred to in paragraph \ref{pLfunc} appears as a special instance of $\mu_E\{x\ri y\}$ (\cite[Lemma 3.10]{HL}).\\

\subsection{Main Conjecture}\label{IMC}
From classical Iwasawa theory, one would expect that 
an analytic $p$-adic $L$-function should annihilate the dual Selmer group $\S$ of the
corresponding extension $\F/F$. 
 
Notations are as in section \ref{pLf}. We let $\wt\F$ be either an anticyclotomic extension as in
\cite[Section 2.3]{Lo} or a cyclotomic one as in \cite[Section 3]{Lo} or in \ref{cycext} above.
The torsion of $G=Gal(\wt\F/H)$ is a finite subgroup and, unless 
possibly when $K$ is a quadratic field where the place $\infty$ splits, $\G:=G/G_{tor}$ is
the Galois group of a $\Zpn$-subextension $\F/H$. Let $\pi:\Z[[G]]\rightarrow\Z[[\G]]$ 
be the natural projection.

\begin{conj}\label{MC}
Assume that $E[p^\infty](\F)=0$ or that $Fitt_{\L_d}(\S_d)$ is principal for all $d$.
For all the $L(E)$'s as described above one has 
\[ \wt{Fitt}_\L(\S)=(\pi(L(E)))\ .\]
\end{conj} 

\begin{rem}
{\em 
1. The interesting cases are when the Fitting ideal is not zero. By analogy with 
the classical situation one expects that then the group $E(\F)$ should be finitely generated. 
The only case we are aware of where the behaviour of Mordell-Weil has been studied in a 
$\Zpn$-tower of function fields is \cite{Br}: Breuer proves that $E(H[\p^{\infty}])$ has 
infinite rank. Observe that no analytic $L$-function has been defined in his setting.\\
2. The available (partial) proofs of Main Conjecures over number fields often require to modify 
local conditions defining Selmer groups. We made a first attempt in Definition \ref{DefModSel}
in order to get torsion $\L_d$-modules. To tackle Conjecture \ref{MC} one might have to further
refine local conditions at places of bad reduction.\\  
3. The analytic $p$-adic $L$-function is an element 
$L(E)\in\Z[[\G]]\,$. In the prospect of Iwasawa theory, it seems natural to 
extend scalars to $\Zp$; however, $L(E)$ may be seen as belonging to $\Zl[[\G]]$ 
as well and one can wonder if it also annihilates the Selmer dual arising from the $l$-torsion. 
In absence of any support or clue, this is just a speculation.}
\end{rem}

\noindent A. Bandini\\
Universit\`a della Calabria - Dipartimento di Matematica\\
via P. Bucci - Cubo 30B - 87036 Arcavacata di Rende (CS) - Italy\\
bandini@mat.unical.it

\noindent I. Longhi\\ 
Universit\`a degli Studi di Milano - Dipartimento di Matematica\\
via Cesare Saldini 50 - 20133 Milano MI - Italy\\
longhi@mat.unimi.it

\end{document}